\documentclass[11pt,leqno]{amsart}

\usepackage[T1]{fontenc}
\usepackage[utf8]{inputenc}
\usepackage{lmodern}
\usepackage[margin=27mm]{geometry}

\usepackage{amsmath,amssymb,amsthm,mathtools}
\usepackage{mathrsfs}
\usepackage{enumitem}
\usepackage{microtype}
\usepackage{xcolor}

\usepackage[
colorlinks=true,
linkcolor=blue!45!black,
citecolor=blue!45!black,
urlcolor=blue!45!black,
pdfstartview=FitV
]{hyperref}

\allowdisplaybreaks[1]
\numberwithin{equation}{section}

\setlist[enumerate]{
	leftmargin=*,
	itemsep=3pt,
	topsep=5pt
}
\setlist[itemize]{
	leftmargin=*,
	itemsep=3pt,
	topsep=5pt
}

\theoremstyle{plain}

\newtheorem{thrm}{Theorem}

\newtheorem{theorem}{Theorem}[section]

\newtheorem{lemma}[theorem]{Lemma}
\newtheorem{proposition}[theorem]{Proposition}
\newtheorem{conjecture}[theorem]{Conjecture}

\theoremstyle{definition}

\newtheorem{example}[theorem]{Example}

\theoremstyle{remark}
\newtheorem{remark}[theorem]{Remark}

\DeclareMathOperator{\GL}{GL}
\DeclareMathOperator{\Hom}{Hom}

\DeclareMathOperator{\Irr}{Irr}

\DeclareMathOperator{\Rep}{Rep}

\DeclareMathOperator{\St}{St}

\newcommand{\one}{\mathbf{1}}

\newcommand{\Z}{\mathbb{Z}}
\newcommand{\cl}{\mathcal{C}_{\mathrm{L}}}

\newcommand{\onto}{\twoheadrightarrow}
\newcommand{\into}{\hookrightarrow}

\newcommand{\seg}[2]{[#1,#2]_2}
\newcommand{\sg}[1]{[#1]_2}

\title[Classification of
$\mathrm{GL}_1(D)\times\mathrm{GL}_{n-1}(D)$-distinguished
representations]
{Classification of
	$\mathrm{GL}_1(D)\times\mathrm{GL}_{n-1}(D)$-distinguished
	representations}

\author[Quanyong Chen]{Quanyong Chen}
\address{
	The University of Hong Kong, China}
\email{chenquanyong@hrbeu.edu.cn}
\author[Jiaqi Huang]{Jiaqi Huang}
\address{ Harbin Engineering University, Harbin, China
}
\email{jiaqihuang@hrbeu.edu.cn}
\begin{document}
	
	\begin{abstract}
		Let $D$ be a quaternion division algebra over a non-Archimedean
		local field $F$ of characteristic zero. 
		We prove the conjecture of Dagar, Sharma, and Verma classifying irreducible smooth
		representations of $\mathrm{GL}_n(D)$ distinguished by
		$\mathrm{GL}_1(D)\times\mathrm{GL}_{n-1}(D)$ for $n>2$.
		These are precisely the trivial representation and the
		irreducible normalized parabolic inductions
		$\one_{n-2}\times\tau$, where $\tau$ is an infinite-dimensional
		irreducible representation of $\mathrm{GL}_2(D)$ distinguished
		by $\mathrm{GL}_1(D)\times\mathrm{GL}_1(D)$.
		The proof combines a choice of cuspidal inducing factor
		with orbit conditions and multisegment methods.
	\end{abstract}
	
	\maketitle
	
	\section{Introduction}
	
	Distinction problems study the existence of invariant linear
	forms on representations under the action of a subgroup.
	For reductive groups and their symmetric subgroups, these
	questions are closely related to harmonic analysis on
	homogeneous spaces. Linear periods provide an important
	example of this connection. In this paper, we classify
	representations of general linear groups over a quaternion
	division algebra admitting a linear period with respect
	to a block-diagonal subgroup of type $(1,n-1)$.
	
	Let $F$ be a non-Archimedean local field of characteristic zero,
	and let $D$ be a quaternion division algebra over $F$.
	For $n\geq 1$, set $G_n=\mathrm{GL}_n(D)$.
	Given nonnegative integers $p,q$ with $p+q=n$, let
	\begin{align*}
		H_{p,q}
		=
		\left\{
		\operatorname{diag}(g_1,g_2):
		g_1\in G_p,\ g_2\in G_q
		\right\}
		\subseteq G_n.
	\end{align*}
	An irreducible smooth complex representation $\pi$ of $G_n$
	is called $H_{p,q}$-distinguished if ${\rm Hom}_{H_{p,q}}(\pi,\one)\neq 0$.
	The elements of this space are called linear periods.
	Our concern is the subgroup $H_n=H_{1,n-1}$.
	
	The theory of linear periods involves both uniqueness
	and classification questions.
	For split general linear groups, Jacquet and Rallis \cite{JR96} established the uniqueness of linear periods for block-diagonal subgroups.
	Chen and Sun \cite{CS20} studied twisted linear periods for the pair
	$(\mathrm{GL}_{2n}(F),
	\mathrm{GL}_n(F)\times\mathrm{GL}_n(F))$
	and their relation to twisted Shalika periods.
	On the classification side, Yang \cite{Yang22}
	characterized the irreducible unitary representations
	of $\mathrm{GL}_{2n}(F)$ admitting linear periods
	with respect to $\mathrm{GL}_n(F)\times\mathrm{GL}_n(F)$.
	These results illustrate the interaction between
	distinction, parabolic induction, and the structure
	of the irreducible dual.
	
	For quaternionic general linear groups, the rank-two
	problem was studied by Raghuram \cite{Raghuram}
	and Gan--Takeda \cite{GanTakeda}.
	Their results describe the representations of $G_2$
	admitting an $H_{1,1}$-period.
	In higher ranks, Anandavardhanan, Lu, Matringe,
	S\'echerre, and Yang \cite{ALMSY24} established results
	on multiplicities and vanishing of linear periods.
	In particular, their results exclude
	$H_n$-distinction for essentially square-integrable
	representations of $G_n$ when $n>2$.
	This makes the analysis of parabolically induced
	representations central to the classification problem.
	
	Dagar, Sharma, and Verma \cite{DSV} proposed that
	$H_n$-distinction in every rank $n>2$ is governed
	by distinguished representations of $G_2$.
	They proved their conjecture for $n=3$ and $n=4$.
	We recall its statement, writing $\one_r$ for
	the trivial representation of $G_r$ and $\times$
	for normalized parabolic induction.
	
	\begin{conjecture}[{\cite[Conjecture 1.1]{DSV}}]
		\label{conj:main}
		Let $n>2$, and let $\pi\in\Irr(G_n)$.
		Then $\pi$ is $H_n$-distinguished if and only if either
		$\pi\simeq\one_n$ or
		\[
		\pi\simeq\one_{n-2}\times\tau,
		\]
		where $\tau$ is an infinite-dimensional irreducible
		$H_{1,1}$-distinguished representation of $G_2$.
	\end{conjecture}
	
	The induction in the second alternative is necessarily
	irreducible, since it is isomorphic to $\pi$.
	Thus the conjecture describes every nontrivial irreducible
	$H_n$-distinguished representation by adjoining a trivial
	factor to distinguished data of rank two.
	
	The main result of this paper establishes the conjecture
	for all $n>2$.
	
	\begin{thrm}\label{thm:main}
		Conjecture~\ref{conj:main} holds for every $n>2$.
	\end{thrm}
	
	\subsection*{Outline of the proof}
	
	Our proof combines the orbit conditions of \cite{DSV}
	with the structure of representations described by
	multisegments.
	The relevant classification originates in the work
	of Zelevinsky \cite{Ze80} for split general linear groups.
	For general linear groups over division algebras,
	we use the results of Tadi\'c \cite{Tadic}
	and M\'inguez--S\'echerre \cite{MS13}.
	We also use the M\oe glin--Waldspurger algorithm
	\cite{MW86}, extended to $\mathrm{GL}_n(D)$ by
	Badulescu and Renard \cite{BR07}, to pass between
	the Zelevinsky and Langlands parametrizations.
	
	A key step is the choice of a cuspidal inducing factor.
	For $n\geq 5$, we prove that every irreducible non-cuspidal
	representation $\pi$, after replacement by $\pi^\vee$
	if necessary, admits a quotient presentation
	$\rho\times\eta\twoheadrightarrow\pi$,
	where $\rho$ is cuspidal, $\eta$ is irreducible, and
	$\rho\not\simeq\nu^{n-3}$ as a representation of $G_1$.
	This follows by comparing the endpoints of the segments
	parametrizing $\pi$ and $\pi^\vee$.
	The resulting choice excludes one closed-orbit condition
	and permits a uniform argument in higher ranks.
	
	Suppose that $\pi$ is distinguished.
	Its period pulls back to $\rho\times\eta$, so the orbit
	conditions apply to this induced representation.
	If $\rho$ has rank at least two, cuspidality excludes
	the open orbit, while the remaining closed-orbit condition
	and the vanishing results of \cite{ALMSY24} force its rank
	to be exactly two.
	This leads to the representations predicted by
	Conjecture~\ref{conj:main}.
	
	For a rank-one factor, the analysis reduces to the open
	orbit. Frobenius reciprocity then gives an embedding
	$\eta\hookrightarrow\rho^\vee\times\one_{n-2}$,
	so that $\pi$ is an irreducible constituent of $I_n(\rho)=\rho\times\rho^\vee\times\one_{n-2}$.
	We determine the possible distinguished constituents
	of this family by separating the boundary values
	$\rho=\nu^{\pm(n-1)}$ from the remaining cases.
	At the boundary, the M\oe glin--Waldspurger algorithm
	and a vanishing statement for inductions involving
	$\nu^{n-2}\St_2$ exclude every nontrivial constituent.
	Away from the boundary, linked-segment arguments
	and the irreducibility of unitary parabolic induction
	\cite{Secherre} reduce the problem to distinction in
	rank two.
	The converse follows from the sufficient closed-orbit
	condition of \cite{DSV}.
	
	\subsection*{Organization of the paper}
	
	Section~2 fixes the notation and recalls the necessary
	background on representations and linear periods.
	Section~3 constructs quotient presentations with suitable
	cuspidal first factors.
	Section~4 analyzes the distinguished constituents of
	the induced representations arising from the orbit
	conditions, including the boundary cases.
	Section~5 combines these results to complete the proof
	of Theorem~\ref{thm:main}.
	\subsection*{Acknowledgements}
	Quanyong Chen expresses his gratitude to Kei Yuen Chan for valuable discussions and comments. 
	The authors acknowledge the assistance of ChatGPT in preparing this manuscript. The authors have reviewed and verified the AI-assisted content and take full responsibility for the
	mathematical arguments, results.

	\section{Conventions and preliminaries}	
	Fix a non-Archimedean local field $F$ of characteristic zero, and let $D$ be a quaternion division algebra over $F$, so that ${\rm dim}_FD=4$. 
	For each $n \in \mathbb N$, write
	\begin{align*}
		G_n={\rm GL}_n(D).
	\end{align*}
	Let ${\rm Nrd}_{D/F}: G_n \rightarrow  F^\times$
	denote the reduced norm map, and define
	\begin{align*}
		\nu_n(g)=|{\rm Nrd}_{D/F}(g)|_F, \ \ g \in G_n.
	\end{align*}
	For $p,q\geq 0$ with $p+q=n$, let $H_{p,q}$ be the
	block-diagonal subgroup
	\begin{align*}
		H_{p,q}=\left\{\begin{pmatrix}g_1&0\\0&g_2\end{pmatrix}:g_1\in G_p,\ g_2\in G_q\right\} \subseteq G_n.
	\end{align*}
	This subgroup is the fixed-point group of the involution
	\begin{align*}
		\theta(g)=\delta_{p,q}g\delta_{p,q}^{-1},
		\qquad
		\delta_{p,q}=\operatorname{diag}(I_p,-I_q).
	\end{align*}
	
	We equip $H_{p,q}$ with the family of characters
	$\chi_s$ defined by
	\begin{align*}
		\chi_s\bigl(\operatorname{diag}(g_1,g_2)\bigr)=\nu(g_1)^{2s}\nu(g_2)^{-2s}, \qquad s \in \mathbb R.
	\end{align*}
	Equivalently, if $\alpha_s:F^\times\to\mathbb C^\times$
	is given by $\alpha_s(x)=|x|_F^{2s}$, then $\chi_s$ is
	the composition of $\alpha_s$ with the homomorphism
	\begin{align*}
		H_{p,q}\longrightarrow F^\times,
		\qquad
		\operatorname{diag}(g_1,g_2)
		\longmapsto
		{\rm Nrd}(g_1){\rm Nrd}(g_2)^{-1}.
	\end{align*}
	\subsection{Representations of $G_n$}
	In this subsection, as in \cite{BZ}, we review the requisite facts on the representation theory of ${G}_n$. 
	Let ${\Rep}({G}_n)$ denote the category of smooth complex representations of ${G}_n$ of finite length.
	We denote by ${\Irr}({G}_n)$ the set of equivalence classes of irreducible representations, and by $\mathfrak{C}(G_n)$ the subset consisting of cuspidal representations. 
	In this paper, all representations are smooth. 
	For $\pi \in {\Rep}(G_n)$, we denote its contragredient by $\widetilde{\pi}$.
	
	In this paper, whenever the value of $n$ is clear from the context, we write ${G}$ for ${G}_n$.
	Let $P$ be a parabolic subgroup of ${G}_n$ with Levi decomposition $P=MN$, where $M$ is a Levi subgroup and $N$ the unipotent radical.
	Let $\rho$ be a representation of $M$. By inflation, we regard 
	$\rho$ as a representation of $P$, with the unipotent radical 
	$N$ acting trivially. 
	The normalized parabolic induction of $\rho$ to $G$ is defined by
	${\rm Ind}_P^G(\delta_P^{1/2}\otimes \rho)$, where $\delta_P$ denotes the modular character of $P$.
	Let $\alpha=(n_1,\cdots, n_r)$ be a partition of $n$, and let $M_\alpha\cong G_{n_1} \times \cdots \times G_{n_r}$ be the corresponding Levi subgroup of $G_n$.
	Every representation of $M_\alpha$ can be written in the form
	\begin{align*}
		\rho_1\otimes \cdots \otimes \rho_r,  \ \  \rho \in \Rep(G_{n_i}).
	\end{align*}
	We use the notation $\rho_1\times \cdots \times \rho_r$ for the normalized parabolic induction ${\rm Ind}_{P_\alpha}^{G_n}(\rho_1\otimes \cdots \otimes \rho_r)$.
	Parabolic induction is exact and preserves finite length.
	Permuting inducing factors preserves their Jordan--H\"older constituents,
	but does not in general preserve which constituent is a subrepresentation
	or a quotient.
	
	We next recall the Jacquet functor. Let $(\pi, V)$ be a representation of $G$, and let $N$ be the unipotent radical of 
	$P$.  Define $V(N)={\rm Span}\{ \pi(n)v-v: n \in N, v \in V\}$, and set $V_N=V/V(N)$. 
	The group $P$ acts $V_N$ by $$\pi_n(p)(v+V(N))=\delta_P^{-1/2}(p)\pi(p)v+V(N),  \ \ p \in P.$$
	Since $\delta_P$ is trivial on $N$, this action factors through the quotient $P/N$. 
	Hence $V_N$ may be viewed as a representation of the Levi subgroup
	$M\cong P/N$.
	The resulting representation $(\pi_N, V_N)$ is called the normalized Jacquet module of $\pi$.
	For the Levi subgroup $M=G_l \times G_{n-l}$, 
	we denote the corresponding Jacquet module by $r_{l,n-l}(\pi)$.
	
	\subsection{Segments and multisegments}
	
	For a cuspidal representation $\rho$ of $G_d$, write $\nu_\rho=\nu^{s(\rho)}$
	for its reducibility character.
	For $a\leq b$ with $b-a\in\Z_{\geq 0}$, a segment is a set
	\begin{align*}
		[a,b]_{(\rho)}
		=
		\{\nu_\rho^a\rho,\nu_\rho^{a+1}\rho,\ldots,\nu_\rho^b\rho\}.
	\end{align*}
	Its degree is $d(b-a+1)$.
	
	Following \cite[Section 2]{DSV}, we denote by $Z(\Delta)$
	the unique irreducible subrepresentation of the ascending
	cuspidal induction attached to $\Delta$, and by $L(\Delta)$
	its unique irreducible quotient.
	On the trivial cuspidal line of $G_1$, one has
	$\nu_{\one_1}=\nu^2$.
	To distinguish ordinary $\nu$-exponents from
	$\nu_\rho$-exponents, we write
	\[
	\seg{a}{b}
	=
	\{\nu^a,\nu^{a+2},\ldots,\nu^b\},
	\qquad
	\sg{a}=\{\nu^a\},
	\]
	where $b-a\in 2\Z_{\geq 0}$.
	In particular,
	\begin{equation}\label{eq:segment-examples}
		Z(\seg{1-r}{r-1})=\one_r,
		\qquad
		L(\seg{c-1}{c+1})=\nu^c\St_2.
	\end{equation}

	Let $\Delta_1=[a_1, b_1]_{(\rho_1)}$ and $\Delta_2=[a_2,b_2]_{(\rho_2)}$
	be two segments. We say that $\Delta_1$ precedes $\Delta_2$ if,
	after concatenating the terms of $\Delta_1$ and $\Delta_2$ in this
	order, one can select a subsequence which forms a segment of length
	strictly larger than both $\ell(\Delta_1)$ and $\ell(\Delta_2)$.
	The segments $\Delta_1$ and $\Delta_2$ are called linked if
	$\Delta_1$ precedes $\Delta_2$ or $\Delta_2$ precedes $\Delta_1$.

	Let $\mathcal S$ denote the set of segments and $\mathcal M(\mathcal S)$ the set of finite multisegments. 
	Let $\mathfrak m=\{\Delta_1,\ldots,\Delta_s\}$ be a multisegment, and define
	$\lambda(\mathfrak a)=Z(\Delta_1)\times\cdots\times Z(\Delta_s)$, where $\Delta_i$ does not precede $\Delta_j$ for $i<j$.
	Then $\lambda(\mathfrak m)$ has a unique irreducible
	subrepresentation $Z(\Delta_1,\ldots,\Delta_s)$.
	This representation depends only on the multisegment $\mathfrak m$.
	Similarly, set $\pi(\mathfrak a)=L(\Delta_1)\times\cdots\times L(\Delta_s)$.
	The unique irreducible quotient of $\pi(\mathfrak a)$ is denoted by $L(\mathfrak m)$.
	
	We shall use the following consequences of the multisegment
	classification; see \cite[Section 2]{DSV} and
	\cite{BR07, MS13, Tadic}.
	The following lemma provides a realization of
	$Z(\mathfrak m)$ as the unique irreducible quotient of a suitably
	ordered parabolic induction.
	\begin{lemma}\cite[Theorem 2.4]{EA16} \label{quotientclassfification}
		Let $\mathfrak m=\{\Delta_1,\ldots,\Delta_t\}$ be a multisegment, ordered so that $\Delta_j$ does not precede $\Delta_i$ whenever $i<j$.
		Then the normalized parabolic induction
		\begin{align*}
			Z(\Delta_1)\times\cdots\times Z(\Delta_t)   
		\end{align*}
		has $Z(\mathfrak m)$ as its unique irreducible quotient.
		In particular, for $\Delta=[a,b]_{(\rho)}$, there is asurjection
		\begin{align*}
			\nu_\rho^b\rho\times\nu_\rho^{b-1}\rho
			\times\cdots\times\nu_\rho^a\rho\onto Z(\Delta).
		\end{align*}
	\end{lemma}
	We will also use the following compatibility between the
	multisegment parametrization and contragredients.
	\begin{lemma}\cite[Proposition 5.10]{MS13}\label{prop:contragredient}
		For every multisegment $\mathfrak m=\{\Delta_1,\ldots,\Delta_r\}$,
		\begin{align*}
			Z(\mathfrak m)^\vee\simeq Z(\mathfrak m^\vee),
		\end{align*}
		where
		\begin{align*}
			\mathfrak m^\vee=\{\Delta_1^\vee,\ldots,\Delta_r^\vee\}, \qquad [a,b]_{(\rho)}^\vee=
			[-b,-a]_{(\rho^\vee)}.
		\end{align*}
		In particular, $Z(\Delta)^\vee\simeq Z(\Delta^\vee)$.
	\end{lemma}
	We next recall the standard exact sequences associated with a pair of
	linked segments. 
	\begin{lemma}\cite[Lemma 5.12]{MS13}\label{prop:linked-ses}
		Let $\Delta\prec\Delta'$ be linked segments and set $U=\Delta\cup\Delta', V=\Delta\cap\Delta'$. 
		If $V=\varnothing$, we agree that $Z(V)$ is the trivial representation of $G_0$, which may be omitted in inductions. 
		Then there are short exact sequences
		\begin{equation}\label{exacta}
			0\longrightarrow Z(U)\times Z(V)\longrightarrow Z(\Delta)\times Z(\Delta')\longrightarrow Z(\Delta,\Delta')\longrightarrow0 
		\end{equation}
		and
		\begin{align}\label{exactb}
			0\longrightarrow Z(\Delta,\Delta')\longrightarrow Z(\Delta')\times Z(\Delta)\longrightarrow Z(U)\times Z(V)\longrightarrow 0. 
		\end{align}
	\end{lemma}
	For $\mathfrak m, \mathfrak  n \in\mathcal M(\mathcal S)$, write $\mathfrak m \leq \mathfrak a$ if there is a finite chain $\mathfrak m=\mathfrak m_0, \mathfrak m_1,\cdots, \mathfrak a_r=\mathfrak n$	such that each $\mathfrak m_{i+1}$ is obtained from $\mathfrak m_i$ by replacing two linked segments by their union
	and intersection, with an empty omitted.
	Such a replacement is called an elementary operation.
	This relation defines a partial order on $\mathcal M(\mathcal S)$.	
	We recall the following properties of the multisegment
	classification from \cite{MS13}.
	\begin{lemma}\label{prop:composition-factors}
		Let $\mathfrak m=\{\Delta_1,\ldots,\Delta_r\}$, and set
		\begin{align*}
			I(\mathfrak m)=Z(\Delta_1)\times\cdots\times Z(\Delta_r),
		\end{align*}
		where the factors may be taken in any order.
		For a multisegment $\mathfrak n$, denote by
		$[I(\mathfrak m):Z(\mathfrak n)]$ the multiplicity of
		$Z(\mathfrak n)$ as a composition factor of $I(\mathfrak m)$.
		Then the following statements hold:
		\begin{enumerate}[label=(\roman*)]
			\item
			$[I(\mathfrak m):Z(\mathfrak n)]>0$ if and only if $\mathfrak n\leq\mathfrak m$.
			\item
			$[I(\mathfrak m):Z(\mathfrak m)]=1$.
			\item
			If the segments of $\mathfrak m$ are pairwise unlinked, then $I(\mathfrak m)\simeq Z(\mathfrak m)$ is irreducible .
			\item
			If $\pi_1=Z(\mathfrak m)$ and $\pi_2=Z(\mathfrak n)$ are supported on two distinct cuspidal lines, then
			\begin{align*}
				\pi_1\times\pi_2\simeq Z(\mathfrak m+\mathfrak n)
			\end{align*}
			is irreducible.
		\end{enumerate}
	\end{lemma}
	\subsection{M\oe glin--Waldspurger algorithm} 
	The M\oe glin--Waldspurger algorithm provides a combinatorial
	description of the correspondence between the Zelevinsky and
	Langlands parametrizations. Originally established in
	\cite{MW86}, it was extended to $\GL_n(D)$ by Badulescu and
	Renard \cite{BR07}. An additional account of the algorithm
	and its applications to Speh representations can be found
	in \cite{Ba14}.
	The algorithm associates to each multisegment $\mathfrak m$
	a multisegment $\mathfrak m^t$ satisfying
	$Z(\mathfrak m)\simeq L(\mathfrak m^t)$ and $(\mathfrak m^t)^t=\mathfrak m$.
	We recall the recursive construction below.
	
	First suppose that all segments of $\mathfrak m$ lie on
	a single cuspidal line, with endpoints expressed in powers
	of $\nu_\rho$. 
	For a segment $\Delta=[a,b]_{(\rho)}$, set $\Delta^-=[a,b-1]_{(\rho)}$
	interpreting $\Delta^-$ as empty when $a=b$.
	Let $d$ be the largest right endpoint occurring in
	$\mathfrak m$. Among the segments ending at $d$, choose
	$\Delta_0$ with the largest left endpoint. 
	Having chosen $\Delta_{s-1}$, choose $\Delta_s$ with the largest left
	endpoint among those segments that end at $d-s$ and
	precede $\Delta_{s-1}$. Stop when no such segment exists.
	If the resulting chain is
	$\Delta_0,\Delta_1,\ldots,\Delta_r$,
	record the segment $\Gamma=[d-r,d]_{(\rho)}$.
	
	Form a new multisegment $\mathfrak m'$ by replacing each
	selected occurrence of $\Delta_s$ by $\Delta_s^-$,
	discarding empty segments and leaving all other occurrences
	unchanged. The construction is then repeated on
	$\mathfrak m'$. Equivalently,
	$\mathfrak m^t=\{\Gamma\}+(\mathfrak m')^t$, 
	where $+$ denotes segments addition.
	The total length decreases at each step, so the procedure
	terminates. 
	For a multisegment supported on several cuspidal
	lines, the construction is carried out separately on each
	line, and the resulting multisegments are combined.

	\subsection{Finite-length arguments}
	In this subsection, we record two elementary facts about finite-length representations.
	These allow us to extract irreducible inducing factors while
	retaining a given irreducible quotient, and concerns the irreducible
	quotients of a representation with a unique irreducible quotient
	up to isomorphism.
	
	\begin{proposition}\label{lem:tail}
		Let $\rho$ be irreducible, let $X$ have finite length, and let $\theta$ be irreducible. 
		If $\rho\times X\onto\theta$, then there is an irreducible constituent $\eta$ of $X$ such that $\rho\times\eta\onto\theta$.
	\end{proposition}
	\begin{proof}
		Choose a composition series $0=X_0\subset X_1\subset\cdots\subset X_t=X$.
		For a fixed surjection $f:\rho\times X\to\theta$, take the least $i$ such that $f$ is nonzero on $\rho\times X_i$.
		That restriction is surjective because $\theta$ is irreducible, and it vanishes on $\rho\times X_{i-1}$.
		Exactness of induction gives the required surjection from $\rho\times(X_i/X_{i-1})$.
	\end{proof}

	\begin{proposition}\label{lem:head}
		Suppose a finite-length representation $M$ has a unique irreducible quotient up to isomorphism, equal to $L$.
		Every nonzero quotient of $M$ has an irreducible quotient isomorphic to $L$.   
	\end{proposition}
	\begin{proof}
		Every nonzero finite-length representation has an irreducible quotient.
		Any irreducible quotient of a quotient of $M$ is itself a quotient of $M$.
	\end{proof}

	\subsection{Results about distinction}
	The external period results needed here are the following.
	\begin{enumerate}[label=(P\arabic*)]
		\item\label{item:contra}
		If an irreducible $\pi$ is $H_n$-distinguished, so is $\pi^\vee$ \cite[Lemma 2.6]{DSV}.
		\item\label{item:DSvan}
		An essentially square-integrable representation of $G_r, r>2$, has no $H_{1,r-1}$-period against any character \cite[Theorem 3.8]{ALMSY24}.
		\item\label{item:ranktwo}
		For an infinite-dimensional $\tau\in {\rm Irr}(G_2)$,
		\begin{equation}\label{eq:rank-two}
			{\rm Hom}_{H_{1,1}}(\tau,\chi_s)\ne0 \quad\Longleftrightarrow\quad
			{\rm Hom}_{H_{1,1}}(\tau,\one)\ne0 \qquad s\in\mathbb R.
		\end{equation}
		We use the formulation of \cite[Remark 3.1]{DSV}. The rank-two period theory is developed in \cite{Raghuram,GanTakeda}.
		\item\label{item:unitary}
		Normalized induction of an irreducible unitary representation of a Levi subgroup of $G_r$ is irreducible \cite{Secherre}.
	\end{enumerate}
	
	Let $\pi_1\in {\rm Irr}(G_k)$ and $\pi_2\in {\rm Irr}(G_{n-k})$.
	The necessary orbit conditions of \cite[Lemma 2.8]{DSV} say that
	$\pi_1\times\pi_2$ can be $H_n$-distinguished only if at least
	one of the following holds:
	\begin{enumerate}[label=(\alph*)]
		\item\label{orb:a}
		$\pi_1=\nu^{n-k-2}\one_k$, and $\pi_2$ is
		$(H_{1,n-k-1},\chi_{k/2})$-distinguished.
		\item\label{orb:b}
		$\pi_1$ is $(H_{1,k-1},\chi_{(k-n)/2})$-distinguished, and
		$\pi_2=\nu^{2-k}\one_{n-k}$.
		\item \label{orb:c}
		Define $L_{n,k}=G_{k-1}\times\Delta(G_1)\times G_{n-k-1}$,
		Then one has
		\begin{equation}\label{eq:open-general}
			\Hom_{L_{n,k}}\left(
			r_{(k-1,1)}(\pi_1)\otimes r_{(1,n-k-1)}(\pi_2),
			\nu^{n-k-1}\boxtimes\one \boxtimes\nu^{-(k-1)}
			\right)\ne0,
		\end{equation}
	\end{enumerate}
	where $\Delta(G_1)=\{(g,g): g\in G_1\} \subset G_1\times G_1$.
	The two closed conditions also give sufficient conditions for distinction \cite[Lemma 2.9]{DSV}.

	\begin{proposition}\label{lem:open}
		Let $m=n-2, \sigma \in {\rm Irr}(G_1)$, and $\eta\in {\rm Irr}(G_{m+1})$.
		If the open orbit contributes to \(\sigma\times\eta\), then
		\begin{align*}
			\eta\hookrightarrow\sigma^\vee\times\one.
		\end{align*}
	\end{proposition}	
	\begin{proof}
		Denote by $S$ and $V$ the spaces of $\sigma$ and $\eta$, respectively, and set
		\begin{align*}
			J=r_{(1,m)}(\eta).
		\end{align*}
		Denote the action of $G_1\times G_m$ on $J$ by $j$.
		For a rank-one first factor, the open-orbit condition gives a nonzero linear functional
		\begin{align*}
			\ell:S\otimes J\longrightarrow\mathbb C
		\end{align*}
		such that
		\begin{align*}
			\ell\bigl(\sigma(a)s\otimes j(a,B)x\bigr)=\ell(s\otimes x)
		\end{align*}
		for all $a\in G_1, B\in G_m, s\in S$, and $x\in J$.
		
		Every irreducible smooth representation of \(D^\times\) is finite-dimensional, so its smooth dual equals its algebraic dual.
		Define
		\begin{align*}
			T:J\longrightarrow S^\vee, \qquad T(x)(s)=\ell(s\otimes x).
		\end{align*}
		The invariance of \(\ell\) implies
		\begin{align*}
			T(j(a,B)x)(s)=\ell(s\otimes j(a,B)x)=\ell(\sigma(a^{-1})s\otimes x)=(\sigma^\vee(a)T(x))(s).
		\end{align*}
		Thus $T$ is a homomorphism
		\begin{align*}
			T:J\longrightarrow\sigma^\vee\boxtimes\mathbf1_m.
		\end{align*}
		Moreover, $T\neq0$ since $T=0$ would imply $\ell=0$.
		
		We now make the Frobenius reciprocity map explicit.
		Let $G=G_{m+1}$ and let
		\begin{align*}
			P=P_{(1,m)}=
			\left\{
			\begin{pmatrix}a&y\\0&B\end{pmatrix}: a\in D^\times,\ B\in G_m,\ y\in M_{1,m}(D)
			\right\}.
		\end{align*}
		Write $P=MN$, where $M\simeq G_1\times G_m$ is the block-diagonal Levi subgroup and $N$ is the unipotent radical.
		Let $q:V\to J$ be the natural quotient map.
		By the normalization of the Jacquet module,
		\begin{align*}
			q(\eta(p)v)=\delta_P(p)^{1/2}j(a,B)q(v), \qquad p=\begin{pmatrix}a&y\\0&B\end{pmatrix}.
		\end{align*}	
		Define
		\begin{align*}
			\Phi(v)(g)=T\bigl(q(\eta(g)v)\bigr), \qquad v\in V,\quad g\in G.
		\end{align*}
		Then
		\begin{align*}
			\Phi(v)(pg)=\delta_P(p)^{1/2}\sigma^\vee(a)\Phi(v)(g),
		\end{align*}
		which is precisely the transformation law for the normalized induction $\sigma^\vee\times\mathbf1_m$.
		If an open compact subgroup fixes $v$, it also fixes $\Phi(v)$ under right translation, so $\Phi(v)$ is smooth.
		
		Furthermore,
		\begin{align*}
			\Phi(\eta(g_0)v)(g)=\Phi(v)(gg_0),
		\end{align*}
		hence $\Phi$ is a $G$-homomorphism.
		It is nonzero. Indeed, since $q$ is surjective and $T\neq0$, some $v\in V$ satisfies
		\begin{align*}
			\Phi(v)(1)=T(q(v))\neq0.
		\end{align*}
		Since \(\Phi\neq0\) and $\eta$ is an irreducible representation, we have $\eta\hookrightarrow\sigma^\vee\times\mathbf1_m$.
	\end{proof}
	\section{Distinguished constituents of induced representations}
	
	In this section, we study the distinguished irreducible
	constituents of $\rho\times\rho^\vee\times\one_{N-2}$.
	A vanishing result for inductions involving a Steinberg twist
	allows us to treat the boundary cases. 
	
	\subsection{Choosing the first cuspidal factor}
	For an irreducible non-supercuspidal representation $\pi$ of $G_n$,
	define
	\begin{align*}
		\cl(\pi)=
		\left\{\rho\in {\rm Cusp}(G_r):
		1\le r<n,\ 
		\rho\times\eta\onto\pi
		\text{ for some }\eta\in \Irr(G_{n-r})\right\}.
	\end{align*}
	
	\begin{proposition}\label{lem:endpoint}
		Let $\pi=Z(\mathfrak m)\in {\rm Irr}(G_n)$ be non-cuspidal.
		On any cuspidal line occurring in $\mathfrak m$, choose a segment $\Delta=[a,b]_{(\rho)}$ with minimal ending.
		Then $\nu_\rho^b\rho\in\cl(\pi)$. 
	\end{proposition}
	
	\begin{proof}
		Let $\rho\in \mathfrak{C}(G_d)$.
		Any segment preceding $\Delta$ must lie on the same cuspidal
		line and have a strictly smaller ending. 
		Thus no segment of $\mathfrak m$ precedes $\Delta$.
		
		We may therefore order the segments as
		$\Delta_1=\Delta,\Delta_2,\ldots,\Delta_t$, with endings nondecreasing on each cuspidal line. 
		No later segment precedes an earlier one. 
		By Lemma \ref{quotientclassfification}, the quotient form of the multisegment classification gives
		\begin{align*}
			Z(\Delta)\times Z(\Delta_2)\times\cdots\times Z(\Delta_t)
			\twoheadrightarrow\pi.
		\end{align*}
		Using the surjection
		\begin{align*}
			\nu_\rho^b\rho\times\cdots\times\nu_\rho^a\rho
			\twoheadrightarrow Z(\Delta)
		\end{align*}
		and exactness of induction, we obtain
		\begin{align*}
			f:\rho_0\times X\twoheadrightarrow\pi,
			\qquad \rho_0=\nu_\rho^b\rho,
		\end{align*}
		where
		\begin{align*}
			X=
			\nu_\rho^{b-1}\rho\times\cdots\times\nu_\rho^a\rho
			\times Z(\Delta_2)\times\cdots\times Z(\Delta_t),
		\end{align*}
		with empty lists omitted. This is a finite-length representation
		of $G_{n-d}$, since normalized parabolic induction preserves
		finite length.
		
		Choose a composition series
		\begin{align*}
			0=X_0\subset X_1\subset\cdots\subset X_r=X,
		\end{align*}
		and let $i$ be minimal such that
		$f(\rho_0\times X_i)\neq0$.
		Irreducibility of $\pi$ implies that this image is all of
		$\pi$, while minimality gives
		$f(\rho_0\times X_{i-1})=0$.
		Exactness of induction therefore yields a surjection
		\begin{align*}
			\rho_0\times(X_i/X_{i-1})\twoheadrightarrow\pi.
		\end{align*}
		The constituent $X_i/X_{i-1}$ is irreducible.
		
		Finally, $d<n$. otherwise $\mathfrak m$ would consist of a
		single segment of length one, making $\pi$ cuspidal.
		Since $\rho_0$ is cuspidal, the displayed surjection
		proves $\rho_0\in\cl(\pi)$.
	\end{proof}

	\begin{proposition}\label{lem:choice}
		Let $n\ge5$, and let $\pi\in{\rm Irr}(G_n)$ be non-cuspidal.
		For at least one of $\pi$ and $\pi^\vee$, there is a presentation
		\begin{align*}
			\rho\times\eta\onto\pi
			\quad\text{or}\quad
			\rho\times\eta\onto\pi^\vee
		\end{align*}
		with $\rho$ cuspidal, $\eta$ irreducible, and
		$\rho\not\simeq\nu^{n-3}$ as a representation of $G_1$.
	\end{proposition}
	\begin{proof}
		We argue by contradiction and assume that the conclusion fails.
		Proposition~\ref{lem:endpoint} shows that
		both sets of possible first factors are nonempty, and hence
		\begin{equation}\label{eq:bad-only}
			\cl(\pi)=\cl(\pi^\vee)=\{\nu^{n-3}\}.
		\end{equation}
		Apply Proposition~\ref{lem:endpoint} separately to each cuspidal line of the
		multisegment $\mathfrak m$ for $\pi$.
		Every such line must contain $\nu^{n-3}$. 
		Thus all segments lie on the single step-two line containing this character.
		
		Set $A=n-3$, and write
		$\mathfrak m=\{\seg{a_i}{b_i}\}_{i=1}^t$in ordinary
		$\nu$-exponents. Minimal endings for $\pi$, followed by minimal
		endings for
		$\pi^\vee=Z(\{\seg{-b_i}{-a_i}\}_{i=1}^t)$, give
		\begin{align*}
			\min_i b_i=A,\qquad -\max_i a_i=A.
		\end{align*}
		Every segment therefore contains $\seg{-A}{A}$, and has degree
		at least $A+1=n-2$. Since $2(n-2)>n$, there is only one segment.
		Its endpoints must then be exactly $-A,A$, giving total degree
		$n-2$, a contradiction.
	\end{proof}
	
	\begin{remark}
		Lemma~\ref{lem:choice} is independent of distinction and of the central character. 
		Its role is to eliminate orbit \ref{orb:a} for a suitable
		cuspidal first factor. This is why the proof below does not require
		an inductive classification of twisted periods in smaller ranks.
	\end{remark}
	
	\subsection{Vanishing and the boundary configuration}
	In this subsection, we establish a vanishing result for
	inductions involving $\nu^{N-2}\St_2$. Combining this with
	the orbit conditions and the M\oe glin--Waldspurger algorithm,
	we show that the trivial representation is the only
	$H_N$-distinguished irreducible constituent of
	$\nu^{-(N-1)}\times\one_{N-2}\times\nu^{N-1}$.
	
	\begin{proposition}\label{lem:steinberg}
		Let $N>2$ and $\eta\in {\rm Irr}(G_{N-2})$. Then
		\begin{align*}
			{\rm Hom}_{H_N}
			\bigl(\nu^{N-2}\operatorname{St}_2\times\eta,\one \bigr)=0.
		\end{align*}
	\end{proposition}
	\begin{proof}
		Let $\tau=\nu^{N-2}\operatorname{St}_2$.
		Suppose that $\tau\times\eta$ is $H_N$-distinguished.
		The three-orbit criterion with $k=2$ then requires at least
		one of the following conditions.
		
		Closed orbit \ref{orb:a} requires
		\begin{align*}
			\tau\simeq\nu^{N-4}\mathbf1_2,
		\end{align*}
		which is impossible because $\tau$ is infinite-dimensional.
		
		Closed orbit \ref{orb:b} requires
		\begin{align*}
			{\rm Hom}_{H_{1,1}}
			\bigl(\tau,\chi_{(2-N)/2}\bigr)\neq0.
		\end{align*}
		However, for $z\in F^\times$,
		\begin{align*}
			\chi_{(2-N)/2}(zI_2)=1,
			\qquad
			\tau(zI_2)=|z|_F^{4(N-2)}{\rm id},
		\end{align*}
		since \(\operatorname{St}_2\) has trivial central character.
		As \(N>2\), these central characters differ, so this Hom
		space vanishes.
		
		For the open orbit, by  \cite[Proposition 4.10]{MS13}, the normalized Jacquet module is
		\begin{align*}
			r_{(1,1)}(\tau)=(\nu^{N-2}\boxtimes\nu^{N-2}) \otimes(\nu\boxtimes\nu^{-1})=\nu^{N-1}\boxtimes\nu^{N-3}.
		\end{align*}
		Consequently, the first $G_1$-factor of
		$G_1\times\Delta G_1\times G_{N-3}$ acts on the source
		by $\nu^{N-1}$, whereas the target character restricts
		to $\nu^{N-3}$. These characters are distinct, so the open-orbit Hom space is trivial.
		
		Thus none of the three necessary conditions holds, we complete the proof.
	\end{proof}
	
	We now apply the preceding vanishing result to determine
	the distinguished irreducible constituents in the following boundary case.
	\begin{proposition}\label{lem:boundary}
		Let $N\ge4$ and $a=-(N-1), b=N-1$, and set
		\begin{align*}
			B_N=\nu^a\times Z(\seg{a+2}{b-2})\times\nu^b
			=\nu^a\times\one_{N-2}\times\nu^b.
		\end{align*}
		The only $H_N$-distinguished irreducible constituent of $B_N$
		is $\one_N$.
	\end{proposition}
	\begin{proof}
		The elementary operations starting from
		$\{\sg a,\seg{a+2}{b-2},\sg b\}$ can attach the left endpoint,
		the right endpoint, or both endpoints to the middle segment.
		The possible constituents are therefore precisely
		\begin{align}
			T_N&=Z(\seg a b)=\one_N,\label{eq:boundary-list}\\
			U_N&=Z(\sg a,\seg{a+2}b),\nonumber\\
			V_N&=Z(\seg a{b-2},\sg b)=U_N^\vee,\nonumber\\
			W_N&=Z(\sg a,\seg{a+2}{b-2},\sg b).\nonumber
		\end{align}
		No multiplicity assertion is needed here.
		
		\smallskip
		\noindent\emph{Excluding $U_N$ and $V_N$.}
		Set $E=Z(\sg a,\seg{a+2}{b-2})$, an irreducible representation
		of $G_{N-1}$. The M\oe glin--Waldspurger algorithm gives
		\begin{align}
			U_N&=L(\sg b,\sg{b-2},\ldots,\sg{a+4},\seg a{a+2}),\label{eq:MW-U}\\
			E&=L(\sg{b-2},\ldots,\sg{a+4},\seg a{a+2}).\label{eq:MW-E}
		\end{align}
		To see the algorithm explicitly, remove the right endpoints of the
		long segment one at a time. 
		Each step produces a singleton until the final remaining endpoint
		is combined with $\sg{a}$ to form $\seg{a}{a+2}$.
		The two computations differ only by the initial singleton $\sg b$.
		
		Let
		\begin{align*}
			R=\nu^{b-2}\times\nu^{b-4}\times\cdots\times\nu^{a+4}
			\times L(\seg a{a+2}).
		\end{align*}
		The factors are in Langlands order. Equations
		\eqref{eq:MW-U}--\eqref{eq:MW-E} say that $R$ has unique
		irreducible quotient $E$, while $\nu^b\times R$ has unique
		irreducible quotient $U_N$. 
		Exactness gives a nonzero quotient
		$\nu^b\times R\onto\nu^b\times E$. Accoring to Proposition~\ref{lem:head}, we have
		\begin{equation}\label{eq:U-presentation}
			\nu^b\times E\onto U_N.
		\end{equation}
		
		If $U_N$ were distinguished, so would be the representation on
		the left of \eqref{eq:U-presentation}. 
		Since the first inducing factor is $\nu^b=\nu^{N-1}$,
		neither closed-orbit condition can hold, as they require this
		factor to be $\nu^{N-3}$ and $\nu^{-(N-1)}$, respectively.
		By Proposition~\ref{lem:open}, the open-orbit embedding would require
		\begin{equation}\label{eq:impossible-E}
			E\into\nu^{-b}\times\one_{N-2}
			=\nu^a\times Z(\seg{a+2}{b-2}).
		\end{equation}
		The induction on the right has length two and unique irreducible
		socle $Z(\seg a{b-2})$, by the linked-segment rule. 
		Its other constituent is $E$, a different irreducible representation.
		Hence \eqref{eq:impossible-E} is impossible.
		This excludes $U_N$, and property \ref{item:contra} excludes
		$V_N=U_N^\vee$.
		
		\smallskip
		\noindent\emph{Excluding $W_N$.}
		The algorithm gives
		\begin{equation}\label{eq:MW-W}
			W_N=L(\seg{b-2}b,\sg{b-4},\sg{b-6},\ldots,
			\sg{a+4},\seg a{a+2}).
		\end{equation}
		The first step joins $\sg b$ to the right endpoint of the middle
		segment. Subsequent steps remove singletons, and the final step joins
		the left endpoint to $\sg a$. When $N=4$, the list of middle
		singletons in \eqref{eq:MW-W} is empty.
		
		Thus $W_N$ is a quotient of an induction
		$L(\seg{b-2}b)\times X$ with $X$ of finite length.
		From Proposition~\ref{lem:tail}, replaces $X$ by an irreducible representation $\eta$.
		But
		\begin{align*}
			L(\seg{b-2}b)=\nu^{b-1}\St_2=\nu^{N-2}\St_2,
		\end{align*}
		so Proposition~\ref{lem:steinberg} excludes distinction of $W_N$.
		Only $T_N=\one_N$ remains.
	\end{proof}
	
	\begin{example}[The boundary formulas in rank five]
		For $N=5$, the boundary cuspidal exponents are
		$-4,-2,0,2,4$, and
		\begin{align*}
			\begin{aligned}
				U_5&=Z(\sg{-4},\seg{-2}4)
				=L(\sg4,\sg2,\seg{-4}{-2}),\\
				E&=Z(\sg{-4},\seg{-2}2)
				=L(\sg2,\seg{-4}{-2}),\\
				W_5&=Z(\sg{-4},\seg{-2}2,\sg4)
				=L(\seg24,\sg0,\seg{-4}{-2}).
			\end{aligned}
		\end{align*}
		In particular $W_5$ is a quotient of
		$\nu^3\St_2\times\one_1\times\nu^{-3}\St_2$, to which
		Lemma~\ref{lem:steinberg} applies after extracting an irreducible tail.
	\end{example}
	
	\subsection{The rank-one open-orbit family}
	We study distinction for $\one_m\times\tau$ in terms of
	$H_{1,1}$-periods on $\tau\in\Irr(G_2)$. Together with
	the boundary case treated above, this criterion allows us
	to determine the distinguished irreducible constituents of
	$\rho\times\rho^\vee\times\one_{N-2}$.
	
	\begin{proposition}\label{lem:character}
		For $m\ge2$ and $\tau\in{\rm Irr}(G_2)$,
		\begin{align*}
			{\rm Hom}_{H_{1,m+1}}(\one_m\times\tau,\one)\ne0
		\end{align*}
		if and only if $\tau$ is infinite dimensional and
		$H_{1,1}$-distinguished. 
	\end{proposition}
	\begin{proof}
		Apply the orbit conditions with \(n=m+2\) and \(k=m\).
		Orbit \ref{orb:b} would require
		\begin{align*}
			\Hom_{H_{1,m-1}}(\one_m,\chi_{-1})\ne0,
		\end{align*}
		which is false.
		For orbit \ref{orb:c}, by the definition of the normalized Jacquet functor \cite[Definition 1.22]{Solleveld}, we have
		\begin{align*}
			r_P(\one_m)
			=(\nu^{-1}\one_{m-1})\boxtimes\nu^{m-1}.
		\end{align*}
		The first factor would have to be $\nu\one_{m-1}$. The characters
		are different since $m-1\ge1$.
		Hence only orbit \ref{orb:a} can contribute, and it requires
		$\Hom_{H_{1,1}}(\tau, \chi_{m/2})\ne0$.
		
		Every finite-dimensional irreducible smooth representation of $G_2$ is a character. Such a character restricts to $H_{1,1}$ with the same character on both factors, so it cannot equal
		$\nu^m\boxtimes\nu^{-m}$. Thus $\tau$ is infinite dimensional,
		and \eqref{eq:rank-two} gives its $H_{1,1}$-distinction.
		Conversely, \eqref{eq:rank-two} and the sufficiency of closed
		orbit \ref{orb:a} give the desired period on $\one_m\times\tau$.
	\end{proof}
	
	The following proposition describes the possible distinguished
	constituents of $I_N(\rho)$.
	\begin{proposition}\label{lem:initial}
		Let $N\ge4$, and let $\rho\in {\rm Irr}(G_1)$.
		Every $H_N$-distinguished irreducible constituent of
		\begin{align*}
			I_N(\rho)=\rho\times\rho^\vee\times\one_{N-2}
		\end{align*}
		is either $\one_N$, or an irreducible induction
		$\one_{N-2}\times\tau$ with $\tau\in\Irr(G_2)$ infinite
		dimensional and $H_{1,1}$-distinguished.
	\end{proposition}
	\begin{proof}
		Set $m=N-2$, and write
		\begin{align*}
			\Delta_m=\seg{-(m-1)}{m-1},
			\qquad Z(\Delta_m)=\one_m.
		\end{align*}
		If $\rho=\nu^{\pm(m+1)}$, invariance of semisimplified
		parabolic induction under permutation of the factors shows that
		$I_N(\rho)$ has the same irreducible constituents as $B_N$.
		The assertion then follows from Proposition~\ref{lem:boundary}.
		
		Assume henceforth that $\rho\ne\nu^{\pm(m+1)}$.
		A singleton is linked to $\Delta_m$ only when it is immediately
		adjacent to one of its endpoints, namely
		$\sg{-(m+1)}$ or $\sg{m+1}$.
		Consequently, neither $[\rho]$ nor $[\rho^\vee]$ is linked to $\Delta_m$.
		The corresponding factors can therefore be interchanged, giving
		\begin{align*}
			I_N(\rho)\simeq\one_m\times Q,
			\qquad Q=\rho\times\rho^\vee.
		\end{align*}
		
		If $Q$ is irreducible, the three segments
		$[\rho],  [\rho^\vee]$, and $\Delta_m$ are pairwise unlinked.
		Thus $I_N(\rho)$ is irreducible, and the assertion follows from Proposition~\ref{lem:character} with $\tau=Q$. 
		
		Suppose now that $Q$ is reducible. The reducibility criterion gives
		\begin{align*}
			\rho^\vee\simeq\nu_\rho^\varepsilon\rho,
			\qquad \varepsilon\in\{1,-1\}.
		\end{align*}
		Set $\sigma=\nu_\rho^{\varepsilon/2}\rho$. Then
		\begin{align*}
			\sigma^\vee
			\simeq\nu_\rho^{-\varepsilon/2}\rho^\vee
			\simeq\nu_\rho^{\varepsilon/2}\rho
			=\sigma.
		\end{align*}
		Twisting does not change the reducibility step, so
		$\nu_\sigma=\nu_\rho$, and therefore
		\begin{align*}
			\{\rho,\rho^\vee\}
			=
			\{\nu_\sigma^{-1/2}\sigma,\nu_\sigma^{1/2}\sigma\}.
		\end{align*}
		Moreover, $Q$ has length two; write its composition series as
		\begin{align*}
			0\longrightarrow\tau_1\longrightarrow Q
			\longrightarrow\tau_2\longrightarrow0,
			\qquad \tau_1,\tau_2\in{\rm Irr}(G_2).
		\end{align*}
		
		If $\sigma\not\simeq\one_1$, the cuspidal line containing
		$\rho$ and $\rho^\vee$ is disjoint from the line supporting
		$\one_m$. 
		Indeed, an intersection would imply
		$\sigma\simeq\nu^t$ for some $t\in\mathbb R$.
		Self-duality would then give $\nu^t=\nu^{-t}$, forcing $t=0$,
		contrary to $\sigma\not\simeq\one_1$.
		Since both $\tau_i$ have cuspidal support
		$\{\rho,\rho^\vee\}$, disjointness of cuspidal lines implies
		that each $\one_m\times\tau_i$ is irreducible.
		
		If $\sigma\simeq\one_1$, then
		$\nu_\sigma=\nu^2$ and $\rho=\nu^{\pm1}$.
		In this case
		\begin{align*}
			\{\tau_1,\tau_2\}=\{\one_2, {\rm St}_2\}.
		\end{align*}
		Both $\one_m\times\one_2$ and
		$\one_m\times {\rm St}_2$ are irreducible by
		the irreducibility of normalized induction of irreducible
		unitary representations.
		
		Thus, in either case, exactness of normalized parabolic induction
		gives
		\begin{align*}
			0\longrightarrow\one_m\times\tau_1
			\longrightarrow I_N(\rho)
			\longrightarrow\one_m\times\tau_2
			\longrightarrow0
		\end{align*}
		with irreducible outer terms. These are therefore precisely the
		irreducible constituents of $I_N(\rho)$.
		Finally, by Proposition~\ref{lem:character}, a constituent
		$\one_m\times\tau_i$ is distinguished only if $\tau_i$ is
		infinite-dimensional and $H_{1,1}$-distinguished.
	\end{proof}

	\section{proof of the conjecture}
	In this section, we combine the preceding results to prove Conjecture~\ref{conj:main}. A suitable choice of cuspidal first factor allows us to apply the orbit conditions and
	reduce necessity to the cases already treated. 
	The sufficiency follows from the closed-orbit criterion.
	\begin{theorem}\label{prop:necessity}
		Let $n>2$, and let $\theta\in {\rm Irr}(G_n)$ be $H_n$-distinguished. Then \(\theta\) belongs to the list in Conjecture~\ref{conj:main}.
	\end{theorem}
	\begin{proof}
		For $n=3,4$, use \cite[Theorems 3.3 and 4.3]{DSV}.
		Assume $n \ge 5 $. Property \ref{item:DSvan} implies that $\theta$ is not cuspidal. 
		By Proposition~\ref{lem:choice}, and after replacing $\theta$ by $\theta^\vee$ if necessary, we have
		\begin{equation}\label{eq:main-quotient}
			\rho\times\eta\onto\theta,\qquad
			\rho\in \mathfrak{C}(G_r),\quad \eta\in {\rm Irr}(G_{n-r}),\quad \rho\not\simeq\nu^{n-3}.
		\end{equation}
		The replacement preserves distinction by \ref{item:contra}.
		The period on \(\theta\) pulls back through \eqref{eq:main-quotient},
		so the orbit conditions apply to \(\rho\times\eta\).
		
		\smallskip
		\noindent\emph{First case: \(r\ge2\).}
		Orbit \ref{orb:a} is impossible because $\rho$ is cuspidal and is not a character. 
		Orbit \ref{orb:c} is impossible because $r_{(r-1,1)}(\rho)=0$. 
		Thus orbit \ref{orb:b} gives
		\begin{equation}\label{eq:rank-r-closed}
			\eta=\nu^{2-r}\one_{n-r},\qquad {\rm Hom}_{H_{1,r-1}}(\rho,\chi_{(r-n)/2})\ne0.
		\end{equation}
		By \ref{item:DSvan}, this forces $r=2$.
		Then $\eta=\one_{n-2}$, and \eqref{eq:rank-two} shows that the cuspidal $\rho$ is $H_{1,1}$ distinguished.
		The induction $\rho\times\one_{n-2}$ is irreducible, since its two factors have disjoint cuspidal lines. 
		Therefore $\theta\simeq\one_{n-2}\times\rho$, as required.
		
		\smallskip
		\noindent\emph{Second case: $r=1$.}
		The choice in \eqref{eq:main-quotient} excludes orbit \ref{orb:a}.
		If orbit \ref{orb:b} held, then
		\begin{align*}
			\rho=\nu^{-(n-1)},\qquad \eta=\nu\one_{n-1}.
		\end{align*}
		In the notation of Proposition~\ref{lem:boundary}, this gives
		\begin{align*}
			\rho\times\eta=\nu^a\times Z(\seg{a+2}b).
		\end{align*}
		Its unique irreducible quotient is $U_n$, by the linked-segment rule.
		Hence $\theta\simeq U_n$, contradicting Proposition~\ref{lem:boundary}. 
		Thus orbit \ref{orb:b} is also impossible.
		
		The open orbit must therefore contribute.
		By	Proposition~\ref{lem:open}, we have $\eta\into\rho^\vee\times\one_{n-2}$, and exact induction gives
		\begin{align*}
			\rho\times\eta\into \rho\times\rho^\vee\times\one_{n-2}.
		\end{align*}
		Since $\theta$ is a quotient of the left-hand side, it is an irreducible constituent of the right-hand side.
		It follows from Proposition~\ref{lem:initial} that $\theta$
		belongs to the required list.
		
		Finally, the list is stable under contragredients.
		The dual of $\one_n$ is $\one_n$, and
		\begin{align*}
			(\one_{n-2}\times\tau)^\vee \simeq\one_{n-2}\times\tau^\vee.
		\end{align*}
		The representation $\tau^\vee$ is again infinite dimensional and $H_{1,1}$-distinguished by \ref{item:contra}.
		Thus the same conclusion holds for the original $\theta$ if we made a contragredient replacement.
		
		We now prove the sufficiency. 
		The trivial representation is distinguished by definition. 
		Suppose $\theta=\one_{n-2}\times\tau$ is irreducible, with $\tau$ infinite dimensional and $H_{1,1}$-distinguished.
		By \eqref{eq:rank-two},
		\begin{align*}
			{\rm Hom}_{H_{1,1}}(\tau,\chi_{(n-2)/2})\ne0.
		\end{align*}
		Apply the sufficient closed-orbit condition \cite[Lemma 2.9(1)]{DSV} with $k=n-2$.
		Its first inducing factor is $\nu^{n-k-2}\one_k=\one_{n-2}$, and its remaining period is
		exactly the one just displayed.
		Hence $\Hom_{H_n}(\theta,\one)\ne0$.
	\end{proof}

	\begin{remark}[Quotient directions used in the argument]
		Throughout the proof, we pull periods back along surjective
		homomorphisms without assuming that they descend from an induced
		representation to each of its irreducible quotients.
	\end{remark}

\end{document}